\documentclass[12pt,reqno]{amsart}
\usepackage{amsmath,amsfonts,amssymb,amscd,amsthm,amsbsy,epsf}
\numberwithin{equation}{section}
\newtheorem{theorem}{Theorem}
\newtheorem{defn}{Definition}

\newtheorem{proposition}[theorem]{Proposition}
\theoremstyle{remark}
\newtheorem*{remark}{Remark}
\newcommand{\reg}{k+\alpha}
\newcommand{\g}{\mathfrak{g}}
\def\Diff{\operatorname{Diff}}
\def\Id{\operatorname{Id}}
\def\Lip{\operatorname{Lip}}

\newcommand{\N}{\mathbb{N}}
\newcommand{\R}{\mathbb{R}}
\def\real{{\mathbb R}}
\def\torus{{\mathbb T}}
\begin{document}
\title[Smooth dependence]
{Smooth dependence on parameters of solution of cohomology equations over 
Anosov systems and applications to cohomology  equations on 
diffeomorphism groups} 
\author{R. de la Llave \and A. Windsor}
\address{Department of Mathematics, The University of Texas at Austin,
Austin, TX 78712} 
\email{llave@math.utexas.edu}
\address{Department of Mathematical Sciences, 
University of Memphis, Memphis, TN 38152} 
\email{awindsor@memphis.edu}
\begin{abstract}
  We consider the dependence on parameters of the solutions of
  cohomology equations over Anosov diffeomorphisms.  We show that the
  solutions depend on parameters as smoothly as the data.  As a
  consequence we prove optimal regularity results for the solutions of
  equations taking value in diffeomorphism groups.  These results are
  motivated by applications to rigidity theory, dynamical systems, and
  geometry.

  In particular, in the context of diffeomorphism groups we show: Let
  $f$ be a transitive Anosov diffeomorphism of a compact manifold $M$.
  Suppose that $\eta \in C^{\reg}(M,\Diff^r(N))$ for a compact
  manifold $N$, $k,r \in \N$, $r \geq 1$, and $0 < \alpha \leq
  \Lip$. We show that if there exists a $\varphi\in
  C^{\reg}(M,\Diff^1(N))$ solving
  \begin{equation*}
    \varphi_{f(x)} = \eta_x \circ \varphi_x
  \end{equation*}
  then in fact $\varphi \in C^{\reg}(M,\Diff^r(N))$.
\end{abstract}

\subjclass[2000]{
58F15,
22E65,
58D05,
37D20,
37C99 
}

\keywords{Cohomology equations, Anosov diffeomorphisms, Liv\v{s}ic theory,
  diffeomorphism groups, rigidity}

\maketitle 
\baselineskip=18pt              
\section{Introduction}

Let $f$ be a diffeomorphism of a compact manifold $M$.  A cohomology
equation over $f$ is an equation of the form
\begin{equation} \varphi_{f(x)} = \eta_x \cdot \varphi_{x}
  \label{eq:cohomology}
\end{equation} where $\eta$ is given and we are to determine
$\varphi$. The interpretation of $\cdot$ varies depending on the
context but is typically either a group operation or the composition
of linear maps. 

Equations of the form \eqref{eq:cohomology} have been studied for many
different maps (including rotations, and horocycle flows).  In this
paper we will always consider $f$ to be an Anosov diffeomorphism,
usually a transitive Anosov diffeomorphism.

It is relevant to point out that if $f^n (p)=p$ then, by applying
\eqref{eq:cohomology} repeatedly, we obtain
\begin{equation*} 
  \varphi_{f^n(p)} = \eta_{f^{n-1}(p)} \cdots \eta_{p}
  \cdot \varphi_{p}.
\end{equation*} 
Hence, if there is a solution $\varphi$ then, since $\varphi_{f^n(p)}
= \varphi_{p}$, we must have
\begin{equation} \label{obstruction}
  \eta_{f^{n-1}p}\cdots \eta_{p} = \Id.
\end{equation}

When $f$ is an Anososv diffeomorphism a great deal of effort, starting
with the pioneering work \cite{Livsic71}, has been devoted to showing
that \eqref{obstruction}, supplemented with regularity and
``localization'' assumptions, is sufficient for the existence and
regularity of $\varphi$.

The main goal of this paper is to study the parameter dependence of
these solutions given that the data depends smoothly on parameters.

In this paper, we will not be concerned with the existence of
solutions of \eqref{eq:cohomology}, rather,  we will assume that a
solution exists and establish smoothness with respect to parameters.
Of course, there is an extensive literature establishing the existence of
solution, see for example the references in \cite{LlaveW07}.

There are several contexts in which to study cohomological equations
\eqref{eq:cohomology}.  The most classical one is when $\varphi,\eta$
take values in a Lie group $G$. In dynamical systems and geometry,
\eqref{eq:cohomology} also appears in some different contexts. Given a
bundle $E$ over $M$, we take $\varphi$ to be an object defined on the
fiber $E_x$ and $\eta_x$ is an action transporting objects in $E_x$ to
objects in $E_{f(x)}$. For example, in \cite{LlaveW07}, one can find
an application where $\varphi_x$ are conformal structures on $E_x$ (a
subspace of $T_xM$) and $\eta_x$ is the natural transportation of the
conformal structure by the differential. We note that it is not
necessary to assume that the bundle $E$ is finite dimensional.  It
suffices to assume that the linear operators lie on a Banach
algebra~\cite{BercoviciN}.  Of course, when the bundle $E$ is trivial,
the linear operators on $E_x$ can be identified with a matrix group,
so that this geometric framework reduces to the Lie group framework
with $G$ a matrix group (or, more generally, a Banach algebra of
operators).

One very interesting example, which indeed serves as the main
motivation for this paper is when $\varphi_x$ and $\eta_x$ are
supposed to be $C^r$ diffeomorphisms of a compact manifold $N$ and the
operation in the right hand side of \eqref{eq:cohomology} is just
composition. Though $\Diff^r (N)$ is certainly a group under
composition, it is not a Lie group since the group operation is not
differentiable as a map from $\Diff^r(N)$ to $\Diff^r(N)$. The problem
was considered in \cite{NiticaT96} where it was shown that results on
\eqref{eq:cohomology} have implications for rigidity. In
\cite{NiticaT96} one can find results on existence of solutions when
$N= \torus^d$ and, in \cite{LlaveW07} for a general $N$.

Both in \cite{NiticaT96}, \cite{LlaveW07}, from the fact that
$\eta_x\in \Diff^r(N)$ one can only reach the conclusion that
$\varphi_x \in \Diff^{r-R} (N)$ (in \cite{LlaveW07} the regularity
loss $R=3$, but in \cite{NiticaT96} $R$ depends on the dimension of
$N$).

The main goal of this paper is to overcome this loss of regularity and
show that if $\eta_x \in \Diff^r(N)$ for $r \geq 1$, and $\varphi \in
\Diff^1(N)$, then $\varphi \in \Diff^r(N)$.  See Theorem~\ref{diffeo}
for a more precise formulation.

As it turns out, the main technical tool in the proof of
Theorem~\ref{diffeo} is to establish results on the smooth dependence on
parameters for the solutions of \eqref{eq:cohomology} which may be of
independent interest.  We formulate them as Theorem~\ref{bootstrapbundle}, 
Theorem~\ref{bootstrapgroup}. 

When $G$ is a commutative group smooth dependence on parameters is
rather elementary, see Proposition~\ref{livsic}, and our main
technique is to reduce to the commutative case.


\subsection{Sketch of the proofs}

In this section we informally present the main ideas behind the proof
omitting several of the details which we will treat later. 

\subsubsection{Some remarks on the relation between bootstrap of
regularity and dependence on parameters}

Going over the proofs in \cite{NiticaT96} and, more explicitly, in
\cite{LlaveW07} it becomes apparent that, the loss of regularity of
the solutions is closely related to the fact that the group operation
is not differentiable.

Hence, we find it convenient to turn the tables.  Rather than
considering $\varphi :M\times N\to N$ as a function from $M$ to a
space of mappings on $N$, we consider $\varphi$ as a mapping from $M$
to $N$ with parameters from $N$.  In this application, $N$ plays two
roles: as the ambient space for the diffeomorphisms and as the
parameter. For the sake of clarity, we will develop our main technical
results denoting the target space as a $N$ and the set of parameters
as $U$. This is, of course, more general, and in many cases, the
parameters appear independently of the target space. 

More precisely, if we fix $ u \in U$ and suppose that $\varphi$ is
differentiable with respect to $u$ then the chain rule tells us that
\begin{equation}
  D_u\varphi_{f(x)} (u)
  = D_u\eta_x \circ \varphi_{x} (u) \cdot D_u\varphi_{x} (u)
  \label{derivative}
\end{equation} 
where $D_u$ denotes the derivative with respect to the parameter $u
\in U$. Note that for typographical convenience we will often use
$\eta_x^u$ in place of $\eta(x,u)$, and $\varphi_x^u$ in place of
$\varphi(x,u)$.

We see that for each $u \in U$ \eqref{derivative} is an equation of
the form \eqref{eq:cohomology} for the cocycle $D_u \varphi_x(u)$ with
generator
\begin{equation}
  \label{bootstrap}
   \tilde\eta_x^u = D_u\eta_x \cdot \varphi_{x}^{u}
\end{equation}
In this case $E$ is the bundle of linear maps from $T_uN$ to
$T_{\varphi_x(u)} N$.

If $\varphi_{x}^{u}$ depends $C^1$ on $u$ and $\eta_x$ depends $C^2$
on $u$, then $\tilde \eta_x$ is $C^1$ in $u$. Using the regularity
theory for commutative cohomology equations we obtain that $D_u
\varphi_{x}^{u}$ is $C^1$ with respect to $u$. This means that
$\varphi_x^u$ is $C^2$ with respect to $u$. This argument to improve
the regularity can be repeated so long as we can differentiate
$\eta_x$.

\subsubsection{The dependence on parameters}

The idea of the proof of the smooth dependence on parameters for
\eqref{eq:cohomology} is more involved that the bootstrap of
regularity since we need to justify the existence of the first
derivative.

For simplicity of notation, we interpret \eqref{eq:cohomology} in the
linear bundle maps framework. This will be the crucial technical
result for the bootstrap of regularity in diffeomorphism groups. The
case of Lie group valued cocycles will be dealt with in
Section~\ref{Liegroups}.

We want to show that if $\varphi^u,\eta^u$ solve
\eqref{eq:cohomology} and $\eta^u$ depends smoothly on
parameters, then $\varphi^u$ depends smoothly on parameters.

However solutions of \eqref{eq:cohomology} are not unique.  Taking
advantage of this it is very easy to construct solutions which depend
badly on parameters.  It is therefore necessary to impose some
additional condition such as that $\varphi_p^u$ is smooth in $u \in U$
for some fixed $p \in M$. For convenience we will take $p \in M$ to be a
periodic point.

In order to prove differentiability we first find a candidate for the
derivative using Liv\v{s}ic methods and then argue that this candidate
is the true derivative. The key observation is that if we formally
differentiate \eqref{eq:cohomology} with respect to the parameter $u$
we obtain
\begin{align} 
  D_u\varphi_{f(x)}^u & = D_u \eta_x^u \cdot \varphi_{x}^u + \eta_x^u
  \cdot D_u\varphi_{x}^u \nonumber
  \intertext{Using \eqref{eq:cohomology} we get} 
  D_u\varphi_{f(x)}^u & = D_u \eta_{x}^{u} \cdot \varphi_{x}^{u} +
  \varphi_{f(x)}^u \cdot (\varphi_x^u)^{-1} \cdot D_u \varphi_{x}^u 
  \label{firststep}
\end{align}
Multiplying on the right both sides of \eqref{firststep} by
$(\varphi_{f(x)}^u)^{-1}$, we obtain
\begin{equation} 
  (\varphi_{f(x)}^u)^{-1} \cdot D_u\varphi_{f(x)}^{u} =
  (\varphi_{f(x)}^{u})^{-1} \cdot D_{u}\eta_{x}^{u} \cdot
  \varphi_{x}^{u} + (\varphi_{x}^u)^{-1} \cdot D_u\varphi_{x}^{u}
\label{derivative2}
\end{equation}
where $(\varphi_{x}^u)^{-1}$ refers to the inverse of the linear map
$\varphi_{x}^u$. This equation \eqref{derivative2} is a cohomology
equation over a commutative group for the new unknown
\begin{equation*}
  \xi_x =  (\varphi_{x}^{u})^{-1} D_u\varphi_{x}^{u}.
\end{equation*}
We will show that the periodic obstruction \eqref{obstruction}
corresponding to \eqref{derivative2} is the derivative with respect to
parameters of the periodic obstruction corresponding to
\eqref{eq:cohomology}.

Hence, applying the results on Liv\v{s}ic equations, we obtain a solution
of \eqref{derivative2}. It remains to show that this candidate is the
true derivative. 

As mentioned before, the solutions of \eqref{derivative2} are not
unique, but, under the hypothesis that $\varphi_p (u)=: \gamma(u)$ is
smooth in $u$ for some $p$, which we made at the outset, it is natural
to impose that the solution to \eqref{derivative2} satisfy the
normalization
$$D_u\varphi_p = D_u\gamma.$$

Once we have a candidate for a derivative, we will use a comparison
argument, Proposition~\ref{uniqueness} to show that indeed it is a
derivative.

Of course, once we have established the existence of the first
derivative, we have also shown that it satisfies \eqref{derivative2}.
Hence, to establish the existence of higher derivatives it suffices to
use the (much simpler) theory of dependence of parameters in
commutative cohomology equations.

The above argument uses very heavily that we are dealing with the
framework of linear operators on bundles.  The adaptation of the above
argument to general Lie groups requires some adaptations of the
geometry, see Section~\ref{Liegroups}.

\section{Smooth dependence on parameters of cohomology equation}

In this section, we make precise the previous arguments and establish
smooth dependence on parameters for the solutions of
\eqref{eq:cohomology}. In Section~\ref{sec:regularity} we make precise
what we mean by smooth dependence on parameters. In
Section~\ref{bundles} we present the results when
\eqref{eq:cohomology} is interpreted as an equation between linear
bundle maps. This section will be crucial for the case of
diffeomorphism groups. In Section~\ref{Liegroups} we present the
results for cocycles taking values in a Lie group.

\subsection{Notation on regularity}
\label{sec:regularity}

We will understand that a function is $C^r$ when it admits $r$
continuous derivatives.  When the spaces we consider are not compact,
we will assume that all the derivatives are bounded.

For $0<\alpha \le \Lip$ we will define real valued $C^\alpha$
functions on a metric space $N$ in the usual way
\begin{equation*} H_\alpha (N) = \sup_{x\ne y} \frac{|f(x) -
f(y)|}{d(x,y)^\alpha}
\end{equation*}

As it is well known, for functions taking values on a manifold, there
are several equivalent definitions of H\"older functions but not a
natural norm \cite{HirschP69}.  For our purposes, any of these
definitions will be enough.

We will now introduce the function spaces we will consider. Both are
of the ``rectangular'' type common in dependence on parameters
arguments~\cite{CabreFL03b}. 

\begin{defn}
  Let $M$ and $N$ be compact smooth manifolds, and $U \subseteq \R^n$
  be an open ball. Let $ 0< \alpha \leq \Lip$ and $k, r \in \N$. Let
  $h: M \times U \rightarrow N$.

  We will write $h \in C^{k+\alpha, r}$ if
  \begin{enumerate}
  \item for all $0 \leq i \leq r$ and all $0 \leq j \leq k$ the
    derivative $D_x^jD_u^i \varphi$ exists and is continuous on $M
    \times U$. 

  \item for all $0 \leq i \leq r$ and all $u \in U$ the derivative
    $D_u^i \varphi( \cdot, u) \in C^\alpha(M)$ and the H\"older
    constant does not depend on $u$. 

  \end{enumerate}
 
  We also define $\hat{C}^{\reg, r}(U, M) = C^r\bigl(U,
  C^{\reg}(M)\bigr)$. That, is the set of functions $h: M \times U
  \rightarrow N$ such that $u \rightarrow h(\cdot, u)$ is $C^r$ when
  we give $\varphi(\cdot, u) $ the $C^{\reg}$ topology.
\end{defn}


The reason to introduce these spaces is that $C^{\reg, r}$ is natural
when performing geometric constructions involving derivatives.  The
space $\hat{C}^{\reg, r} $ is natural when we consider dependence on
parameters of commutative Liv\v{s}ic equations.  Of course, the two
spaces are closely related.  The relation is formulated in the
following proposition.

\begin{proposition}\label{prop:inclusion}
For any $0 < \alpha' < \alpha  \leq \Lip$, $k, r \in \N$,
we have: 
\begin{equation}
  \label{inclusion}
  \hat{C}^{\reg, r}  \subset C^{\reg, r} \subset \hat{C}^{k +\alpha', r}
\end{equation}
\end{proposition}

The simple example $\varphi(u,x) = (x-u)^{\reg}$ satisfies $\varphi
\in C^{\reg, 0}$, $\varphi \notin {\hat C}^{\reg,0}$.  Integrating
with respect to $u$ one can obtain similar examples for all $r \in
\N$.

\begin{proof}
  The first inclusion of \eqref{inclusion} is obvious. 

  Let $\varphi \in C^{k+\alpha,r}$. Let $\alpha' = \alpha \cdot
  \theta$ for $0< \theta < 1$. We wish to show that $\varphi \in
  C^r(U, C^{k+\alpha'}(M))$. This is equivalent to $\partial_u^i
  \varphi \in C^0(U, C^{k+\alpha'}(M)$ for all multi-indices $i$ with
  $0\leq |i| \leq r$. Fix $u \in U$ and consider a compact $K$ with $u \in K
  \subset U$. Restricted to $M \times K$ the function $\varphi$ and
  all its partial derivatives are uniformly continuous. So we have 
  \begin{equation*}
    \|\partial_x^j \partial_u^i \varphi( \cdot, u) 
    - \partial_x^j \partial_u^i \varphi( \cdot, u') \|_{0}
  \end{equation*}
  is controlled by $d(u,u')$, for all multi-indices $j$ with $0\leq
  |j|\leq k$. It remains to show that for $j = k$ we have
  \begin{equation*}
    \|\partial_x^j \partial_u^i \varphi( \cdot, u)
    - \partial_x^j \partial_u^i \varphi( \cdot, u') \|_{\alpha'}
  \end{equation*}
  is controlled. For any
  multi-index $j$ with $|j| = k$ we have from the Kolmogorov-Hadamard
  inequalities \cite{LlaveO99}
  \begin{multline*}
    \|\partial_x^j \partial_u^i \varphi( \cdot, u)
    - \partial_x^j \partial_u^i \varphi( \cdot, u') \|_{\alpha'}\\
    \leq C \|\partial_x^j \partial_u^i \varphi( \cdot, u)
    - \partial_x^j \partial_u^i \varphi( \cdot, u') \|_{\alpha}^\theta
    \|\partial_x^j \partial_u^i \varphi( \cdot, u)
    - \partial_x^j \partial_u^i \varphi( \cdot, u') \|_{0}^{ 1-\theta}
  \end{multline*}
  By definition of $\varphi \in C^{\reg,r}$ we have that
  $\|\partial_x^j \partial_u^i \varphi( \cdot, u)\|_{\alpha}$ is
  bounded independent of $u \in U$ and hence
  \begin{equation*}
    \|\partial_x^j \partial_u^i \varphi( \cdot, u) 
    - \partial_x^j \partial_u^i \varphi( \cdot, u') \|_{0}
  \end{equation*}
  is controlled since $\partial_x^j \partial_u^i \varphi$ is uniformly
  continuous. Consequently   $\varphi \in
  C^r(U , C^{k+\alpha'}(M))$. Note that we do not claim uniform
  continuity in $u$. 
\end{proof} 

Next, we discuss the regularity properties of the commutative
cohomology equation. The results for $k=0$ appear in \cite{Livsic71,
  Livsic72} and the results for $k>0$ appear in \cite{LlaveMM86}. The following
result will be one of our basic tools later.

\begin{proposition}\label{livsic}
Let $M$ be a compact manifold and $f$ a transitive Anosov 
diffeomorphism on $M$. Let  $p \in M$ be a periodic
point for $f$. 

Consider the commutative cohomology equation  
\begin{equation} \label{commutative}
\begin{split} 
& \varphi^u\circ f (x) - \varphi^u(x)= \eta^u(x) \\ 
& \varphi(p, \cdot) \in C^r(U) 
\end{split} 
\end{equation} 
Assume that for all values of $u$, the periodic orbit obstruction with
$\eta^u$ vanishes.  Then,
\begin{enumerate} 
\item if $\eta \in \hat{C}^{\reg, r}$, then $\varphi \in
  \hat{C}^{\reg, r}$, and

\item if $\eta \in {C}^{\reg, r}$, then $\varphi \in {C}^{\reg, r}$.
  
\end{enumerate} 

\end{proposition}

\begin{proof} 
  The first statement is obvious since we are considering the solution
  of a linear equation and \cite{LlaveMM86} showed that the
  solution operator is continuous from $C^{\reg}$ to $C^{\reg}/\real$.

  The second statement follows because, by the previous argument, we
  have that $\varphi \in {\hat C}^{\reg', r}$.  In particular for $|i|
  \le r$, $\partial_u^i \varphi( \cdot, u) \in C^{\reg'}$ but this
  function satisfies
  \[
  \partial_u^i \varphi^u\circ f  - 
  \partial_u^i \varphi^u   = 
  \partial_u^i \eta^u
  \]
  Since we assumed that $\eta \in C^{\reg, r}$ we have that
  $\partial_u^i \eta( \cdot, u) \in C^{\reg}(M) $.  and, by the results
  of \cite{LlaveMM86}, $\partial_u^i \varphi( \cdot, u) \in
  C^{\reg}(M) $.
\end{proof}

\subsection{Cohomology Equations for Bundle Maps}
\label{bundles}

As we have mentioned before, cohomology equations for bundle maps
appear naturally when we consider the linearization of a dynamical
system or cohomology equations for cocycles taking values in
diffeomorphism groups.

\begin{theorem}\label{bootstrapbundle}
  Let $M$ and $N$ be compact smooth manifolds, $U \subseteq \R^n$
  open, and $E$ a bundle over $N$. Let $ 0< \reg \leq 1$ and $r \in
  \N$ with $r \geq 1$. 

  Let $f : M \rightarrow M$ be a transitive $C^r$ Anosov
  diffeomorphism, $\sigma: U \rightarrow N$ with $\sigma \in C^r(N)$
  and $\tau : M \times U \rightarrow N$ with $\tau \in C^{\reg,r}$
  (resp. $\tau \in \hat{C}^{\reg, r}$).

  Suppose that we have linear maps
  \begin{align*}
    \varphi_x^u &: E_{\sigma(u)} \rightarrow E_{\tau(x,u)},\\
    \eta_x^u &: E_{\tau(x,u)} \rightarrow E_{\tau(f(x),u)}.
  \end{align*}
  such that for all $u \in U$ we have 
  \begin{equation}\label{cohomologyparam}
    \varphi_{f(x)}^u = \eta_x^u \cdot \varphi_x^u. 
  \end{equation}

  If $\eta \in C^{\reg,r}$ (resp. $\eta \in \hat{C}^{\reg, r}$) and there
  exists a periodic point $p \in M$ such that $\varphi(p, \cdot) \in
  C^r(U)$ then $\varphi \in C^{\reg,r}$ (resp. $\varphi \in
  \hat{C}^{\reg, r}$).
\end{theorem}








The most difficult case of Theorem~\ref{bootstrapbundle} is the case
when $r=1$.  The higher differentiability cases follow from this one
rather straightforwardly. The case $r=1$ of
Theorem~\ref{bootstrapbundle} will be the inductive step for the
bootstrap of regularity for cohomology equations on diffeomorphism
groups. 

The proof we present of the $r=1$ case uses that $f$ is transitive but
the subsequent bootstrap argument does not require that $f$ is
transitive.

As indicated before, we will obtain a candidate for a derivative and,
then, we argue it is indeed a true derivative.

Smooth dependence is a local question. For that reason we consider a
local trivialization of the bundle $E$ about the points $\sigma(u)$
and $\tau(x,u)$. The base space portion of the map $\varphi_x^u$ is as
smooth as the minimum smoothness of $\sigma$ and $\tau$. We will
concentrate therefore on the smoothness of the linear operator in the
fibers.

\begin{proof}
  We note that if we could take derivatives in \eqref{cohomologyparam}
  then $D_u\varphi_x^u$, the derivative of $\varphi_x^u$ with respect
  to $u$, would satisfy
  \begin{equation}
    D_u\varphi_{f(x)}^u = D_u\eta_x^u \cdot \varphi_{x}^u + \eta_x^u
    \cdot D_u \varphi_{x}^u
    \label{paramderivative}
  \end{equation} 
  Observing that $\eta_x^u = \varphi_{f(x)}^u \cdot ( \varphi_x^u
  )^{-1}$ we get
  \begin{equation}
    (\varphi_{f(x)}^u)^{-1} \cdot D_u\varphi_{f(x)}^u =
    (\varphi_{f(x)}^u)^{-1} \cdot D_u\eta_x^u \cdot \varphi_{x}^u +
    ( \varphi_x^u )^{-1} \cdot D_u \varphi_{x}^u.
  \end{equation} 
  Writing
  \begin{equation} \xi_x^u := (\varphi_x^u)^{-1} \cdot D_u \varphi_x^u
    \label{transform}
  \end{equation} 
  we see that $\xi_x^u$ would satisfy
  \begin{equation}
    \xi_{f(x)}^u  - \xi_{x}^u = (\varphi_{f(x)}^u)^{-1} \cdot D_u\eta_x^u \cdot
    \varphi_{x}^u 
    \label{Wequation}
  \end{equation}

  Note that \eqref{Wequation} is a commutative cohomology equation
  with the group operation being addition on a vector space.

  In this case, as we will show, the methods of \cite{Livsic72,
    LlaveMM86} to obtain existence and regularity for $\xi$.  Using
  \eqref{transform} we will obtain a candidate for $D_u\varphi_x^u$.

  Since we are assuming that there are solutions of
  \eqref{cohomologyparam} for all $u$ in an open set $U$, we have that
  for any $p$, $f^n(p)=p$
  \begin{equation}\label{paramobstruction} \eta_{f^{n-1}(p)}^u \cdot
    \cdots \cdot
    \eta_{p}^u = \Id.
  \end{equation}
  Differentiating \eqref{paramobstruction} with respect to $u$ we
  obtain
  \begin{equation*}
    \begin{split}
      0 = D_u\eta_{f^{n-1}(p)}^u \cdot \eta_{f^{n-2}(p)}^u \cdots
      \eta_{p}^u &+ \eta_{f^{n-1}(p)}^u \cdot D_u
      \eta_{f^{n-2}(p)}^u \cdot \eta_{f^{n-3}(p)}^u \cdots \eta_{p}^u\\
      &+ \cdots + \eta_{f^{n-1}(p)}^u \cdot \eta_{f^{n-2}(p)}^u\cdots
      \eta_{f(p)}^u \cdot D_u\eta_{p}^u.
    \end{split}
  \end{equation*}
  Using the \eqref{cohomologyparam} to express the products of $\eta$,
  we obtain
  \begin{equation}\label{calculation}
    \begin{split}
      0 = D_u \eta_{f^{n-1}(p)}^u \cdot
      \varphi_{f^{n-1}(p)}^u &\cdot (\varphi_{p}^u)^{-1}\\
      &+ \varphi_p^u \cdot (\varphi_{f^{n-1}(p)}^u)^{-1} \cdot D_u
      \eta_{f^{n-2}(p)}^u \cdot
      \varphi_{f^{n-2}(p)}^u \cdot (\varphi_{p}^u)^{-1}\\
      & + \cdots + \varphi_p^u \cdot (\varphi_{f(p)}^u)^{-1} \cdot
      D_u \eta_{p}^u.
    \end{split}
  \end{equation}
  Multiplying \eqref{calculation} on the left by $(\varphi_p^u)^{-1}$ and
  multiplying by $\varphi_p^u$ on the right we obtain
  \begin{equation}\label{calculation2}
    \begin{split}
      0 =(\varphi_p^u)^{-1} \cdot D_u \eta_{f^{n-1}(p)}^u \cdot
      \varphi_{f^{n-1}(p)}^u + (\varphi_{f^{n-1}(p)}^u)^{-1}& \cdot D_u
      \eta_{f^{n-2}(p)}^u \cdot \varphi_{f^{n-2}(p)}^u\\
      &+ \cdots + (\varphi_{f(p)}^u)^{-1} \cdot D_u \eta_{p}^u \cdot
      \varphi_p^u
    \end{split}
  \end{equation}
  which shows that the periodic orbit obstruction for the existence of a
  solution to the cohomology equation \eqref{Wequation}
  vanishes. Hence we have a solution $\xi_x^u$ to \eqref{Wequation}. 

  Using \eqref{transform} we see that we have a candidate for the
  derivative,
  \begin{equation*}
    D_u \varphi_x^u = \varphi_x^u \cdot \xi_x^u.
  \end{equation*}
  It remains to show that this candidate, which we obtained by a
  formal argument is actually the derivative.

  We start by proving that all the solutions of
  \eqref{cohomologyparam} are differentiable in the stable manifold of
  a periodic point. For the sake of notation, and without any loss of
  generality, we will consider a fixed point $p$.

  Applying \eqref{cohomologyparam} repeatedly we have
  \begin{equation}\label{iterate} 
    \varphi_{f^{n+1}(x)}^u = \eta_{f^n(x)}^u \cdots \eta_x^u \cdot
    \varphi_x^u. 
  \end{equation} 
  For $x\in W_p^s$ we have $d(f^n(x),p) \le C_x \lambda^n$ for some
  $0<\lambda<1$, $n\ge 0$. Since $\eta_x^u$ is H\"older in $x$, and
  $\eta_p^u = \Id $ is continuous, in any smooth local trivialization
  around $p$ we obtain
  \begin{equation*}
    \|\eta_{f^n(x)}^u - \Id \| \le C \lambda^{\alpha n}.
  \end{equation*}

  Therefore we can pass to the limit in \eqref{iterate} and obtain for
  $x\in W_p^s$
  \begin{equation}\label{representation}
    \varphi_p^u = \lim_{n\to\infty}\eta_{f^n(x)}^u \cdots \eta_x^u \cdot
    \varphi_x^u 
  \end{equation}
  where the convergence in \eqref{representation} is uniform in
  bounded sets of $W_p^s$ (in the topology of $W_p^s$).

  We will show that \eqref{representation} can be differentiated with
  respect to $u$. By the product rule
  \begin{equation}\label{sum}
    \begin{split}
      D_u \bigl( \eta_{f^n(x)}^u& \cdots \eta_{x}^u \bigr)\\ & =
      D_u\eta_{f^n(x)}^u \cdot \eta_{f^{n-1}(x)}^u \cdots \eta_{x}^u +
      \eta_{f^{n}(x)}^u \cdot D_u\eta_{f^{n-1}(x)}^u \cdot
      \eta_{f^{n-2}
        (x)}^u \cdots \eta_{x}^u\\
      &\qquad
      + \cdots + \eta_{f^{n}(x)}^u \cdots \eta_{f(x)}^u D_u \eta_{x}^u\\
      & = \varphi_{f^{n+1}(x)}^u \cdot (\varphi_{f^{n+1}(x)}^u)^{-1}
      \cdot D_u \eta_{f^n(x)}^u \cdot \varphi_{f^{n}(x)}^u \cdot
      (\varphi_x^u)^{-1} + \varphi_{f^{n+1}(x)}^u \cdot
      (\varphi_{f^{n}(x)}^u)^{-1} \cdot \\
      & \qquad D_u \eta_{f^{n-1}(x)}^u \cdot \varphi_{f^{n-1}(x)}^u
      \cdot (\varphi_x^u)^{-1} + \cdots + \varphi_{f^{n+1}(x)}^u \cdot
      (\varphi_{f^(x)}^u)^{-1} \cdot D_u \eta_{x}^u \cdot
      \varphi_{x}^u \cdot (\varphi_x^u)^{-1}\\
      &= \varphi_{f^{n+1}(x)}^u \cdot \Bigl( \sum_{i=0}^{n-1}
      (\varphi_{f^{i+1}(x)}^u)^{-1} \cdot D_u \eta_{f^i(x)}^u \cdot
      \varphi_{f^i(x)}^{u} \Bigr) \cdot (\varphi_x^u)^{-1}
    \end{split}
  \end{equation} 
  where we have used again \eqref{cohomologyparam}. Since $p$ is a
  fixed point using \eqref{calculation} we obtain that
  \begin{equation*}
    (\varphi_p^u)^{-1} \cdot D_u \eta_p^u \cdot \varphi_p^u =0.
  \end{equation*}
  Because of the assumed regularity on $\varphi_x$, $D_u \eta_x^u$ and
  the exponentially fast convergence of $f^n(x)$ to $p$, we obtain
  that the general term in \eqref{sum} converges to zero
  exponentially.

  By the Weierstrass $M$-test we obtain that $D_u \bigl(
  \eta_{f^n(x)}^u \cdots \eta_{x}^u\bigr)$ converges uniformly on
  compact subsets of $W^s_p$.  Therefore, the limit is the derivative
  of $\lim_{n\to\infty} \eta_{f^n(x)} \cdots \eta_{x}$ and from that,
  it follows immediately that $\varphi_x^u$ is differentiable with
  respect to $u$ for $x\in W_p^s$.

  Of course, the argument so far allows only to conclude that
  $D_u\varphi_x^u$ is bounded on bounded sets of $W_p^s$. This does
  not show that it is bounded on $M$ since $W_p^s$ is unbounded.

  Nevertheless, we realize that the derivative solves
  \eqref{paramderivative} on $W_p^s$.  Of course, so does the
  candidate for the derivative that we produced earlier. {From} the
  fact that these two functions satisfy the same functional equation
  in $W_p^s$, and that they agree at zero, we will show that they
  agree. This will allow us to conclude that the candidate is indeed
  the derivative in $W_p^s$.  Using that it is a continuous function
  on the whole manifold, there is a standard argument that shows it is
  the derivative everywhere. 

\begin{proposition}\label{uniqueness}   
  Let $p \in M$ be a fixed point of $f$. Let $\xi_1,\xi_2$ be
  continuous functions on $W_p^s$.  Assume that $\xi_1 (p) = \xi_2(p)$
  and that both of them satisfy \eqref{Wequation}. Then $\xi_1(x) =
  \xi_2(x)$ for all $x \in W_p^s$.
\end{proposition}

\begin{proof}
  Note that $A(x) \equiv \xi_1 (x) - \xi_2(x)$ satisfies $A(x) =
  A(f(x))$. For every $x \in W^s_p$ we have $\lim_{n \rightarrow
    \infty} f^n(x) = p$. Since $A$ is continuous, and $A(p) =0$ we see
  that
  \begin{equation*}
    \begin{split}
      A(x) &= A \bigl( f^n(x) \bigr) = \lim_{n \rightarrow \infty}
      A\bigl( f^n(x) \bigr)
      = A\bigl( \lim_{n \rightarrow \infty}f^n(x) \bigr)
      = A(p)\\
      &= 0.
    \end{split}
  \end{equation*}
  Thus $\xi_1 (x) = \xi_2(x)$ for all $x \in W^s_x$.
\end{proof}




We have now shown that the function is differentiable in a leaf of the
foliation which is dense. Furthermore, $D_u\varphi_x^{u}$ extends
continuously -- indeed H\"older continuously -- to the whole
manifold. In this case, an elementary real analysis argument, which we
will present now, shows that $\varphi_x^u$ is everywhere
differentiable.

Since we just want to show that the continuous function is indeed a
derivative in the stable leaf, we can just take a trivialization of
the bundles in a small neighborhood. We will use the same notation for
the objects in the trivialization and the corresponding one in the
bundles. This allows us to subtract objects in different fibers.

If $\psi(s)$ is any smooth path in $U$ contained in the trivializing
neighborhood with $\psi(0) =u_0$, $\psi_1 = u_1$ we have, for $x\in
W_p^s$
\begin{equation}\label{incrementformula} 
  \varphi_x^{u_1} -
  \varphi_x^{u_2} = \int_0^1 D_u\varphi_x^{\psi (s)}\, \psi'(s) \,ds
\end{equation}
Both sides of the equation are H\"older in $x$. Therefore, we can pass
to the limit in \eqref{incrementformula} and conclude that the set of
points $x$ for which \eqref{incrementformula} holds is closed.  But we
had already shown that \eqref{incrementformula} held for $x\in W_p^s$,
which is a dense set since $f$ is transitive. Thus
\eqref{incrementformula} holds for every $x \in M$. This shows
immediately that our candidate is indeed the true derivative and
concludes the proof of the case $r=1$ of
Theorem~\ref{bootstrapbundle}.

It is now relatively simple to obtain higher derivatives. The
auxiliary function $\xi$ satisfies the commutative cohomology equation
\eqref{Wequation}. By assumption $D_u\eta \in C^{\reg,r-1}$
(resp. $D_u\eta \in \hat{C}^{\reg,r-1}$) hence, if $\varphi \in C^{\reg,
  m}$ (resp. $\varphi \in \hat{C}^{\reg, m}$) for $m \leq r-1$ then the
right hand side of \eqref{Wequation} is in $C^{\reg,m}$
(resp. $\hat{C}^{\reg,m}$). Applying Proposition \ref{livsic} we then
obtain $D_u\varphi \in C^{\reg, m}$ (resp. $D_u\varphi \in
\hat{C}^{\reg, m}$) and thus $\varphi \in C^{\reg, m+1}$ (resp. $\varphi
\in C^{\reg, m+1}$). The induction stops at $m = r-1$ when we have
exhausted the regularity of $D_u\eta$. At this point $\varphi \in
C^{\reg,r}$ (resp. $\varphi \in \hat{C}^{\reg,r}$) as required.


\end{proof}
\subsection{Cohomology Equations for  Lie Group Valued Cocycles}
\label{Liegroups}

In this section we consider the dependence on parameters of the
solution to \eqref{eq:cohomology} when $\eta$, and $\varphi$ are
functions taking values in a Lie group $G$.  We denote the Lie algebra
of the Lie group $G$ by $\g$ and we denote the identity in $G$ by $e$.

The proof follows along the same lines as the proof of Theorem
\ref{bootstrapbundle}. We derive a functional equation for a candidate
for the first derivative, show that there is a solution for this
equation, and that the candidate is a true derivative. Finally we use
a bootstrapping argument to get full regularity. The main difference
with the previous section is that we have to deal with the fact that
the group operation is not just the product of linear operations, so
that the derivatives of the functional equation involve the
derivatives of the group operation.

We introduce the notation
\begin{equation}\label{notation}
\begin{split} L_g h & = g\cdot h\\ R_g h & = h\cdot g
\end{split}
\end{equation} where $\cdot$ denotes the group operation.

It follows directly from the definitions of $L,R$ that
\begin{equation}\label{leftright}
\begin{split} L_g \circ L_h & = L_{g\cdot h}\\ R_g \circ R_h & =
R_{h\cdot g}
\end{split}
\end{equation}

It is immediate to show that if $g^u, h^u$ are smooth families
\begin{equation}\label{productrule} D_u (g^u \cdot h^u) = DL_{g^u}
(h^u) D_u h^u + DR_{h^u} (g^u) D_u g^u
\end{equation}

\begin{theorem}\label{bootstrapgroup}
  Let $M$ be a compact manifold, $U \subset \R^d$ open, $G$ be a Lie
  group, $f$ a transitive Anosov diffeomorphism of $M$, and $p \in M$
  a periodic point for $f$. Suppose that $\eta: M \times U \rightarrow
  G$ with $\eta \in C^{\reg,r}$ (resp. $\eta \in \hat{C}^{\reg,r}$) for $0 <
  \alpha \leq \Lip$ and $k, r \in \N$ with $r \geq 1$.

  If $\varphi: M \times U \rightarrow G$ solves 
  \begin{equation*}
    \varphi_{f(x)}^u = \eta_x^u \cdot \varphi_x^u
  \end{equation*}
  and $\varphi_p \in C^r(U,G)$ then $\varphi \in C^{\reg, r}$
  (resp. $\varphi \in \hat{C}^{\reg, r}$).
\end{theorem}

The definition of the spaces $C^{\reg,r}$ and $\hat{C}^{\reg,r}$ appears in
Section~\ref{sec:regularity}. 

\begin{proof}
  Taking derivatives of \eqref{eq:cohomology} and applying the product
  rule we obtain that if there is a derivative of $D_u\varphi_x^u$, it
  should satisfy:
  \begin{equation}\label{step1} D_u \varphi_{f(x)}^u = DL_{\eta_x^u}
    (\varphi_x^u) D_u \varphi_x^u + DR_{\varphi_x^u}
    (\eta_x^u) D_u \eta_x^u
  \end{equation}

  We introduce a function $\xi: M \times U \rightarrow \g$ by
  \begin{equation}\label{Wintroduced} D_u \varphi_x^u =
    DL_{\varphi_x^u}
    (e) \xi_x^u
  \end{equation}

  Introducing the notation \eqref{Wintroduced} is geometrically
  natural because we want to transport all the infinitesimal
  derivatives to the identity, so that $\xi_x^u$ takes values in the
  Lie algebra~$\g$.

  In terms of $\xi$, the equation \eqref{step1} becomes
  \begin{equation}\label{step2} DL_{\varphi_{f(x)}^u} (e) \xi_{f(x)}^u
    =
    DL_{\eta_x^u} (\varphi_x^u) DL_{\varphi_x^u} (e) \xi_x^u
    + DR_{\varphi_x^u} (\eta_x^u) D_u \eta_x^u.
  \end{equation}
  The first factor of the first term in \eqref{step2} can be
  simplified
  \begin{equation}\label{step3}
    DL_{\eta_x^u} ( \varphi_x^u) \cdot DL_{\varphi_x^u}
    (e)  = D_u\bigl( L_{\eta_x^u} \circ L_{\varphi_x^u} \bigr) (e)\\
    = D_u\bigl( L_{\eta_x^u \cdot \varphi_x^u}\bigr) (e)\\ 
    = DL_{\varphi_{f(x)}^u} (e).
  \end{equation}
  Substituting \eqref{step3} into \eqref{step2} we obtain
  \begin{equation}\label{step4}
    DL_{\varphi_{f(x)}^u}(e)\, \xi_{f(x)}^u =
    DL_{\varphi_{f(x)}^u}(e) \, \xi_x^u +DR_{\varphi_x^u} (\eta_x^u) D_u
    \eta_x^u.
  \end{equation}
  Multiplying \eqref{step4} in the left by
  $\bigl(DL_{\varphi_{f(x)}^u}(e)\bigr)^{-1}$ we obtain:
  \begin{equation}\label{step5}
    \xi_{f(x)}^u = \xi_x^u +\bigl(DL_{\varphi_{f(x)}^u}(e)\bigr)^{-1}
    DR_{\varphi_x^u} (\eta_x^u) D_u \eta_x^u
  \end{equation}
  Equation \eqref{step5} can be simplified further; from
  \eqref{leftright}, we have
  \begin{equation*}
    \Id = L_{(\varphi_{f(x)}^u)^{-1}} \circ L_{\varphi_{f(x)}^u}
  \end{equation*}
  where $\Id$ is the identity map on $G$. Thus, by the chain rule, we
  have
  \begin{equation*}
    \Id = DL_{(\varphi_{f(x)}^u)^{-1}} (\varphi_{f(x)}^u) \, DL_{\varphi_{f(x)}^u} (e)
  \end{equation*}
  where $\Id$ is the identity map on $\g$. Hence
  \begin{equation}
    \label{step6} 
    \bigl(DL_{\varphi_{f(x)}^u}(e)\bigr)^{-1} = DL_{(\varphi_{f(x)}^u)^{-1}} (\varphi_{f(x)}^u).
  \end{equation}
  Therefore, \eqref{step5} can be rewritten as
  \begin{equation}\label{finalform}
    \xi_{f(x)}^u = \xi_x^u + DL_{(\varphi_{f(x)}^u)^{-1}}
    (\varphi_{f(x)}^u) \, DR_{\varphi_x^u}(\eta_x^u) \, D_u \eta_x^u .
  \end{equation}

  We are assuming $\varphi_p^u$ is a $C^r$ function of $u$. Thus,
  taking derivatives, and using \eqref{Wintroduced} we obtain
  \begin{equation}\label{normalization}
    D_u \varphi_p^u = DL_{\varphi_p^u} (e) \xi_p^u
  \end{equation}
  so we see that $\xi_p : U \rightarrow \g$ is $C^r$.

  In summary, if the function $\varphi_x^u$ is differentiable then the
  $\g$ valued function $\xi$ introduced in \eqref{Wintroduced} would
  satisfy \eqref{finalform}. Equation \eqref{finalform} is a
  cohomology equation for functions taking values in $\g$.  It is
  therefore a commutative cohomology equation. Thus a necessary and
  sufficient condition for the existence of a smooth solution is the
  vanishing of the periodic orbit obstruction.

  In the next paragraphs we will show that indeed this periodic orbit
  obstruction is met, so that indeed one can find a $\xi$ solving
  \eqref{finalform}. Since \eqref{eq:cohomology} holds for all $u$, if
  $f^n (p)=p$, we have for all $u$
  \begin{align*}
    e &= \eta_{f^{n-1}(p)}^u \cdots \eta_p^u\\
    &= (\varphi_p^u)^{-1} \cdot \eta_{f^{n-1}(p)}^u \cdots \eta_p^u
    \cdot \varphi_p^u
  \end{align*}
  Differentiating with respect to $u$, we obtain:
  \begin{equation}\label{firstderivative}
    \begin{split}
      0 = DR_{\eta_{f^{n-1}(p)}^u \cdots \eta_p^u \cdot
        \varphi_p^u}\bigl( (\varphi_p^u)^{-1}& \bigr)
      D_u(\varphi_p^u)^{-1}\\
      + DL_{(\varphi_p^u)^{-1}}(\eta_{f^{n-1}(p)}^u &\cdots \eta_p^u
      \cdot \varphi_p^u) DR_{\eta_{f^{n-2}(p)}^u \cdots \eta_p^u
        \cdot \varphi_p^u}(\eta_{f^{n-1}(p)}^u) D_u\eta_{f^{n-1}(p)}^u \\
      + \cdots +& DL_{ (\varphi_p^u)^{-1} \cdot \eta_{f^{n-1}(p)}^u
        \cdots \eta_{f(p)}^u}(\eta_p^u \varphi_p^u)
      DR_{\varphi_p^u}(\eta_p^u)
      D_u\eta_p^u\\
      &+ DL_ { (\varphi_p^u)^{-1} \cdot \eta_{f^{n-1}(p)}^u \cdots
        \eta_{p}^u}(\varphi_p^u) D_u \varphi_p^u.
    \end{split}
  \end{equation}
  Using \eqref{eq:cohomology} to reduce the products, we transform
  \eqref{firstderivative} into
  \begin{equation}\label{firstderivative2}
    \begin{split}
      0 = DR_{\varphi_p^u}\bigl( (\varphi_p^u)^{-1} \bigr)
      D_u(\varphi_p^u)^{-1}&\\
      + DL_{(\varphi_p^u)^{-1}}( \varphi_p^u)
      DR_{\varphi_{f^{n-1}(p)}^u}&(\eta_{f^{n-1}(p)}^u)
      D_u\eta_{f^{n-1}(p)}^u\\
      + DL_{(\varphi_{f^{n-1}(p)}^u)^{-1}}&( \varphi_{f^{n-1}(p)}^u)
      DR_{\varphi_{f^{n-2}(p)}^u}(\eta_{f^{n-2}(p)}^u)
      D_u\eta_{f^{n-2}(p)}^u \\
      + \cdots + & DL_{ (\varphi_{f(p)}^u)^{-1}}(\varphi_{f(p)}^u)
      DR_{\varphi_p^u}(\eta_p^u) D_u\eta_p^u\\
      &+ DL_{(\varphi_p^u)^{-1}}(\varphi_p^u) D_u \varphi_p^u.
    \end{split}
  \end{equation}
  Finally we observe that
  \begin{equation*}
    DR_{\varphi_p^u}\bigl( (\varphi_p^u)^{-1} \bigr)
    D_u(\varphi_p^u)^{-1}+  DL_{(\varphi_p^u)^{-1}}(\varphi_p^u) D_u
    \varphi_p^u = D_u \bigl((\varphi_p^u)^{-1} \cdot \varphi_p^u \bigr)=0
  \end{equation*}
  and hence \eqref{firstderivative2} becomes
  \begin{equation*}
    \begin{split}
    0 = DL_{(\varphi_p^u)^{-1}}( \varphi_p^u)
    DR_{\varphi_{f^{n-1}(p)}^u}&(\eta_{f^{n-1}(p)}^u)
    D_u\eta_{f^{n-1}(p)}^u\\
    + DL_{(\varphi_{f^{n-1}(p)}^u)^{-1}}&( \varphi_{f^{n-1}(p)}^u)
    DR_{\varphi_{f^{n-2}(p)}^u}(\eta_{f^{n-2}(p)}^u)
    D_u\eta_{f^{n-2}(p)}^u \\+ \cdots +&
    DL_{(\varphi_{f^{2}(p)}^u)^{-1}}( \varphi_{f^{2}(p)}^u)
    DR_{\varphi_{f(p)}^u}(\eta_{f^(p)}^u) D_u\eta_{f(p)}^u \\ 
    &+ DL_{
      (\varphi_{f(p)}^u)^{-1}}(\varphi_{f(p)}^u)
    DR_{\varphi_p^u}(\eta_p^u) D_u\eta_p^u
  \end{split}\end{equation*}
  which shows the vanishing of the periodic orbit obstruction for
  \eqref{finalform}.

  The rest of the argument does not need any modification from the
  argument in the previous case and we just refer to the previous
  section. We argue that $\varphi$ is differentiable in $W_p^s$ for
  some $p \in M$ and that the candidate for the drivative is indeed
  the true derivative.
\end{proof}

\section{Cohomology equations on diffeomorphism groups}

In this section, we consider \eqref{eq:cohomology} when $\eta$ and
$\varphi$ take values in the diffeomorphism group of a compact
manifold.  As before $f$ is a transitive Anosov diffeomorphism. In
order to emphasize that the parameter space and target space are the
same smooth manifold $N$ and to match the notation in~\cite{LlaveW07}
we use $y$ in place of $u$. 

\begin{theorem}\label{diffeo}
  Let $\eta \in C^{\reg}\bigl( M, \Diff^r(N) \bigr)$ and $p$ is a
  periodic point of $f$. 

  If $\varphi \in C^{\reg}\bigl( M, \Diff^1(N) \bigr)$ solves
  \begin{equation}
    \label{eq:DiffeoCohomology}
    \varphi_{f(x)} = \eta_x \circ \varphi_x
  \end{equation}
  and $\varphi_p \in \Diff^r(N)$ then $\varphi \in C^{\reg}\bigl( M,
  \Diff^r(N) \bigr)$. 
\end{theorem}

Notice that questions of differentiability is entirely local and hence
the global structure of the diffemorphism group, which is quite
complicated, does not enter into the argument. This should be compared
with the existence argument in~\cite{LlaveW07} for $k=0$ which devotes
considerable effort to defining the metric on $\Diff^r(N)$. The local
differential structure on $\Diff^r(N)$ can be found
in~\cite{Banyaga97}. We remark that to show that $\varphi: M
\rightarrow \Diff^r(N)$ is $C^\reg$ it suffices to show that all
partial derivatives in $N$ are $C^\reg$ uniformly. 

\begin{proof} 
  Taking derivatives of \eqref{eq:cohomology} with respect to the
  variable $y$ in the manifold $N$ we obtain
  \begin{equation*}
    D_y\varphi_{f(x)} (y) = D_y\eta_x \circ
    \varphi_x (y) \cdot D_y\varphi_{x} (y).
  \end{equation*}
  This is possible since for all $x$ we have $\varphi_x \in
  \Diff^1(N)$. Now we write $\psi_x^y = D_y \varphi_x(y)$ and
  $\hat\eta_x^y = D_y\eta_x \circ \varphi_x (y)$. The equation is then
  \begin{equation*}
    \psi_{f(x)}^y = \hat\eta_x^y \cdot \psi_x^y
  \end{equation*}
  We have
  \begin{align*}
    \psi_x^y &: T_yN \rightarrow T_{\varphi_x(y)}N\\
    \hat\eta_x^y&: T_{\varphi_x(y)}N \rightarrow T_{\varphi_{f(x)}(y)}N.
  \end{align*}
  Thus we have the set up of Theorem \ref{bootstrapbundle} with
  $\sigma(y)=y$ and $\tau(x,y) = \varphi_x^y$. Clearly $\sigma \in
  C^\infty(N)$.

  If we assume that $\varphi \in C^{\reg}(M, \Diff^m(N))$ then $\tau
  \in C^{\reg,m}$. By hypothesis $\eta \in C^{\reg}(M, \Diff^r(N))$
  thus for $m \leq r-1$ $\hat\eta \in C^{\reg,m}$.  Applying Theorem
  \ref{bootstrapbundle} we obtain that $D_y\varphi \in
  C^{\reg,m}$. Thus $\varphi \in C^{\reg, m+1}$. We proceed by
  induction until $m = r-1$ at which point we have $\varphi \in
  C^{\reg, r}$.

  We have addressed on $\varphi_x$ and not $(\varphi_x)^{-1}$ however
  since we know that $\varphi_x \in \Diff^1(N)$ the inverse function
  theorem shows that $(\varphi_x)^{-1}$ is as smooth as $\varphi_x$
  and depends on parameters with the same smoothness.

\end{proof}

\begin{remark} 
  Note that the main result of \cite{LlaveW07} is that if the periodic
  orbit obstruction is met, $f$ is hyperbolic enough, and $\eta \in
  C^{\alpha} \bigl(M,\Diff^r (N)\bigr)$, $r\ge4$ is close enough to
  the identity, then the cohomology equation
  \eqref{eq:DiffeoCohomology} has a solution $\varphi \in
  C^\alpha\bigl(M, \Diff^1(N)\bigr)$. Hence we can apply Theorem
  \ref{diffeo} to show that in fact $\varphi \in C^\alpha\bigl(M,
  \Diff^r(N)\bigr)$.
\end{remark}

\begin{remark}
  Though the existence result requires both a localization assumption
  and a relation between the H\"older exponent and the hyperbolicity
  of $f$ the bootstrap result requires none of these assumptions. 
\end{remark}


\section*{Acknowledgements}

Alistair Windsor gratefully acknowledges the hospitality of the
University of Texas at Austin. R.L. was partially supported by NSF
grants.


\begin{thebibliography}{dlLMM86}

\bibitem[Ban97]{Banyaga97}
Augustin Banyaga.
\newblock {\em The structure of classical diffeomorphism groups}.
\newblock Kluwer Academic Publishers Group, Dordrecht, 1997.

\bibitem[BN98]{BercoviciN}
Hari Bercovici and Viorel Ni{\c{t}}ic{\u{a}}.
\newblock A {B}anach algebra version of the {L}ivsic theorem.
\newblock {\em Discrete Contin. Dynam. Systems}, 4(3):523--534, 1998.

\bibitem[CFdlL03]{CabreFL03b}
Xavier Cabr{\'e}, Ernest Fontich, and Rafael de~la Llave.
\newblock The parameterization method for invariant manifolds. {II}.
  {R}egularity with respect to parameters.
\newblock {\em Indiana Univ. Math. J.}, 52(2):329--360, 2003.

\bibitem[dlLMM86]{LlaveMM86}
R.~de~la Llave, J.~M. Marco, and R.~Moriy{\'o}n.
\newblock Canonical perturbation theory of {A}nosov systems and regularity
  results for the {L}iv\v sic cohomology equation.
\newblock {\em Ann. of Math. (2)}, 123(3):537--611, 1986.

\bibitem[dlLO99]{LlaveO99}
R.~de~la Llave and R.~Obaya.
\newblock Regularity of the composition operator in spaces of {H}\"older
  functions.
\newblock {\em Discrete Contin. Dynam. Systems}, 5(1):157--184, 1999.

\bibitem[dW07]{LlaveW07}
Rafael {de la Llave} and Alistair Windsor.
\newblock Liv\v{s}ic theorems for non-commutative groups including
  diffeomorphism groups and results on the existence of conformal structures
  for {A}nosov systems, 2007.

\bibitem[HP69]{HirschP69}
Morris~W. Hirsch and Charles~C. Pugh.
\newblock Stable manifolds for hyperbolic sets.
\newblock {\em Bull. Amer. Math. Soc.}, 75:149--152, 1969.

\bibitem[Liv71]{Livsic71}
A.~N. Liv{\v{s}}ic.
\newblock Certain properties of the homology of ${Y}$-systems.
\newblock {\em Mat. Zametki}, 10:555--564, 1971.

\bibitem[Liv72]{Livsic72}
A.~N. Liv{\v{s}}ic.
\newblock Cohomology of dynamical systems.
\newblock {\em Izv. Akad. Nauk SSSR Ser. Mat.}, 36:1296--1320, 1972.

\bibitem[NT96]{NiticaT96}
Viorel Ni{\c{t}}ic{\u{a}} and Andrei T{\"o}r{\"o}k.
\newblock Regularity results for the solutions of the {L}ivsic cohomology
  equation with values in diffeomorphism groups.
\newblock {\em Ergodic Theory Dynam. Systems}, 16(2):325--333, 1996.

\end{thebibliography}

\end{document}